\documentclass[11pt]{amsart}
\usepackage{amsmath}
\usepackage{amsxtra}
\usepackage{amscd}
\usepackage{amsthm}
\usepackage{amsfonts}
\usepackage{amssymb}
\usepackage{eucal}

\newtheorem{lem}[subsection]{Lemma}
\newtheorem{mainlemma}[subsection]{Main Lemma}
\newtheorem{prop}[subsection]{Proposition}

\newtheorem{conj}[subsection]{Conjecture}
\newtheorem{thm}[subsection]{Theorem}

\theoremstyle{definition}

\theoremstyle{remark}

\newcommand{\thmref}[1]{Theorem~\ref{#1}}

\newcommand{\lemref}[1]{Lemma~\ref{#1}}
\newcommand{\mainlemref}[1]{Main Lemma~\ref{#1}}
\newcommand{\propref}[1]{Proposition~\ref{#1}}

\newcommand{\conjref}[1]{Conjecture~\ref{#1}}

\newcommand{\nc}{\newcommand}
\nc{\renc}{\renewcommand}
\nc{\ssec}{\subsection}
\nc{\sssec}{\subsubsection}
\nc{\on}{\operatorname}

\nc\ol{\overline}
\nc\wt{\widetilde}
\nc\tboxtimes{\wt{\boxtimes}}
\nc{\alp}{\alpha}

\nc{\ZZ}{{\mathbb Z}}
\nc{\NN}{{\mathbb N}}
\nc{\OO}{{\mathbb O}}
\renc{\SS}{{\mathbb S}}
\nc{\DD}{{\mathbb D}}
\nc{\GG}{{\mathbb G}}

\nc{\Fq}{{\mathbb F}_q}
\nc{\Fqb}{\ol{{\mathbb F}_q}}
\nc{\Ql}{\ol{{\mathbb Q}_\ell}}
\nc{\id}{\text{id}}
\nc\X{\mathcal X}

\nc{\Hom}{\on{Hom}}
\nc{\Lie}{\on{Lie}}
\nc{\Loc}{\on{Loc}}
\nc{\Pic}{\on{Pic}}
\nc{\Bun}{\on{Bun}}
\nc{\IC}{\on{IC}}
\nc{\Aut}{\on{Aut}}
\nc{\rk}{\on{rk}}
\nc{\Sh}{\on{Sh}}
\nc{\Perv}{\on{Perv}}
\nc{\pos}{{\on{pos}}}
\nc{\Conv}{\on{Conv}}
\nc{\Sph}{\on{Sph}}
\nc{\Sym}{\on{Sym}}
\nc{\BunBb}{\overline{\Bun}_B}
\nc{\Buno}{\overset{o}{\Bun}}
\nc{\BunPb}{{\overline{\Bun}_P}}
\nc{\BunBM}{\overline{\Bun}_{B(M)}}
\nc{\BunPbw}{{\widetilde{\Bun}_P}}
\nc{\BunBP}{\widetilde{\Bun}_{B,P}}
\nc{\GUb}{\overline{G/U}}
\nc{\GUPb}{\overline{G/U(P)}}

\nc{\Hhom}{\underline{\on{Hom}}}
\nc\syminfty{\on{Sym}^{\infty}}
\nc\lal{\ol{\lambda}}
\nc\xl{\ol{x}}
\nc\thl{\ol{\theta}}
\nc\nul{\ol{\nu}}
\nc\mul{\ol{\mu}}
\nc\Sum\Sigma
\nc{\oX}{\overset{o}{X}{}}

\nc{\M}{{\mathcal M}}
\nc{\N}{{\mathcal N}}
\nc{\F}{{\mathcal F}}
\nc{\D}{{\mathcal D}}
\nc{\Q}{{\mathcal Q}}
\nc{\Y}{{\mathcal Y}}
\nc{\G}{{\mathcal G}}
\nc{\E}{{\mathcal E}}
\nc{\CalC}{{\mathcal C}}
\nc\Dh{\widehat{\D}}

\nc{\C}{{\mathcal C}}
\nc{\K}{{\mathcal K}}
\renewcommand{\H}{{\mathcal H}}

\nc{\T}{{\mathcal T}}
\nc{\V}{{\mathcal V}}
\renc{\P}{{\mathcal P}}
\nc{\A}{{\mathcal A}}
\nc{\B}{{\mathcal B}}
\nc{\U}{{\mathcal U}}

\nc{\Gr}{\on{Gr}}

\nc{\frn}{{\check{\mathfrak u}(P)}}
\nc{\p}{\overline{\mathfrak p}}
\nc{\q}{\overline{\mathfrak q}}
\nc\f{{\mathfrak f}}

\nc{\qo}{{\mathfrak q}}
\nc{\po}{{\mathfrak p}}
\nc{\s}{{\mathfrak s}}
\nc\w{\text{w}}

\nc\mathi\iota
\nc\Spec{\on{Spec}}
\nc\Mod{\on{Mod}}
\nc{\tw}{\widetilde{\mathfrak t}}
\nc{\pw}{\widetilde{\mathfrak p}}
\nc{\qw}{\widetilde{\mathfrak q}}
\nc{\jw}{\widetilde j}

\nc{\grb}{\overline{\Gr}}
\nc{\I}{\mathcal I}

\nc{\lambdach}{{\check\lambda}}
\nc{\Lambdach}{{\check\Lambda}{}}
\nc{\much}{{\check\mu}}
\nc{\omegach}{{\check\omega}}
\nc{\nuch}{{\check\nu}}
\nc{\etach}{{\check\eta}}
\nc{\alphach}{{\check\alpha}}
\nc{\betach}{{\check\beta}}
\nc{\rhoch}{{\check\rho}}
\nc{\ch}{{\check h}}

\nc{\Hb}{\overline{\H}}

\nc{\BA}{{\mathbb{A}}}
\nc{\BC}{{\mathbb{C}}}
\nc{\BG}{{\mathbb{G}}}
\nc{\BH}{{\mathbb{H}}}
\nc{\BM}{{\mathbb{M}}}
\nc{\BN}{{\mathbb{N}}}
\nc{\BP}{{\mathbb{P}}}
\nc{\BR}{{\mathbb{R}}}
\nc{\BZ}{{\mathbb{Z}}}
\nc{\BV}{{\mathbb{V}}}
\nc{\BW}{{\mathbb{W}}}
\nc{\BU}{{\mathbb{U}}}
\nc{\BX}{{\mathbb{X}}}
\nc{\BY}{{\mathbb{Y}}}
\nc{\BK}{{\mathbb{K}}}

\nc{\CC}{{\mathcal{C}}}
\nc{\CA}{{\mathcal{A}}}
\nc{\CB}{{\mathcal{B}}}

\nc{\CE}{{\mathcal{E}}}
\nc{\CF}{{\mathcal{F}}}
\nc{\CG}{{\mathcal{G}}}
\nc{\CL}{{\mathcal{L}}}
\nc{\CM}{{\mathcal{M}}}
\nc{\CN}{{\mathcal{N}}}
\nc{\CO}{{\mathcal{O}}}
\nc{\CP}{{\mathcal{P}}}
\nc{\CQ}{{\mathcal{Q}}}
\nc{\CR}{{\mathcal{R}}}
\nc{\CS}{{\mathcal{S}}}
\nc{\CT}{{\mathcal{T}}}
\nc{\CU}{{\mathcal{U}}}
\nc{\CV}{{\mathcal{V}}}
\nc{\CW}{{\mathcal{W}}}
\nc{\CZ}{{\mathcal{Z}}}

\nc{\cM}{{\check{\mathcal M}}{}}
\nc{\csM}{{\check{\mathcal A}}{}}
\nc{\oM}{{\overset{\circ}{\mathcal M}}{}}
\nc{\obM}{{\overset{\circ}{\mathbf M}}{}}
\nc{\oCA}{{\overset{\circ}{\mathcal A}}{}}
\nc{\obA}{{\overset{\circ}{\mathbf A}}{}}
\nc{\ooM}{{\overset{\circ}{M}}{}}
\nc{\osM}{{\overset{\circ}{\mathsf M}}{}}
\nc{\vM}{{\overset{\bullet}{\mathcal M}}{}}
\nc{\nM}{{\underset{\bullet}{\mathcal M}}{}}
\nc{\oD}{{\overset{\circ}{\mathcal D}}{}}
\nc{\obD}{{\overset{\circ}{\mathbf D}}{}}
\nc{\oA}{{\overset{\circ}{\mathbb A}}{}}
\nc{\op}{{\overset{\bullet}{\mathbf p}}{}}
\nc{\cp}{{\overset{\circ}{\mathbf p}}{}}
\nc{\oU}{{\overset{\bullet}{\mathcal U}}{}}
\nc{\oZ}{{\overset{\circ}{\mathcal Z}}{}}
\nc{\ofZ}{{\overset{\circ}{\mathfrak Z}}{}}
\nc{\oF}{{\overset{\circ}{\fF}}}

\nc{\fa}{{\mathfrak{a}}}
\nc{\fb}{{\mathfrak{b}}}
\nc{\fg}{{\mathfrak{g}}}
\nc{\fgl}{{\mathfrak{gl}}}
\nc{\fh}{{\mathfrak{h}}}
\nc{\fj}{{\mathfrak{j}}}
\nc{\fq}{{\mathfrak{q}}}
\nc{\fm}{{\mathfrak{m}}}
\nc{\fn}{{\mathfrak{n}}}
\nc{\fu}{{\mathfrak{u}}}
\nc{\fp}{{\mathfrak{p}}}
\nc{\fr}{{\mathfrak{r}}}
\nc{\fs}{{\mathfrak{s}}}
\nc{\fsl}{{\mathfrak{sl}}}
\nc{\hsl}{{\widehat{\mathfrak{sl}}}}
\nc{\hgl}{{\widehat{\mathfrak{gl}}}}
\nc{\hg}{{\widehat{\mathfrak{g}}}}
\nc{\chg}{{\widehat{\mathfrak{g}}}{}^\vee}
\nc{\hn}{{\widehat{\mathfrak{n}}}}
\nc{\chn}{{\widehat{\mathfrak{n}}}{}^\vee}

\nc{\fA}{{\mathfrak{A}}}
\nc{\fB}{{\mathfrak{B}}}
\nc{\fD}{{\mathfrak{D}}}
\nc{\fE}{{\mathfrak{E}}}
\nc{\fF}{{\mathfrak{F}}}
\nc{\fG}{{\mathfrak{G}}}
\nc{\fK}{{\mathfrak{K}}}
\nc{\fL}{{\mathfrak{L}}}
\nc{\fM}{{\mathfrak{M}}}
\nc{\fN}{{\mathfrak{N}}}
\nc{\fP}{{\mathfrak{P}}}
\nc{\fU}{{\mathfrak{U}}}
\nc{\fV}{{\mathfrak{V}}}
\nc{\fZ}{{\mathfrak{Z}}}

\nc{\bb}{{\mathbf{b}}}
\nc{\bc}{{\mathbf{c}}}
\nc{\be}{{\mathbf{e}}}
\nc{\bj}{{\mathbf{j}}}
\nc{\bn}{{\mathbf{n}}}
\nc{\bp}{{\mathbf{p}}}
\nc{\bq}{{\mathbf{q}}}
\nc{\bu}{{\mathbf{u}}}
\nc{\bv}{{\mathbf{v}}}
\nc{\bx}{{\mathbf{x}}}
\nc{\by}{{\mathbf{y}}}
\nc{\bw}{{\mathbf{w}}}
\nc{\bA}{{\mathbf{A}}}
\nc{\bB}{{\mathbf{B}}}
\nc{\bC}{{\mathbf{C}}}
\nc{\bD}{{\mathbf{D}}}
\nc{\bG}{{\mathbf{G}}}
\nc{\bH}{{\mathbf{H}}}
\nc{\bL}{{\mathbf{L}}}
\nc{\bY}{{\mathbf{Y}}}

\nc{\bM}{{\mathbf{M}}}
\nc{\bN}{{\mathbf{N}}}
\nc{\bV}{{\mathbf{V}}}
\nc{\bW}{{\mathbf{W}}}
\nc{\bX}{{\mathbf{X}}}
\nc{\bK}{{\mathbf{K}}}
\nc{\bZ}{{\mathbf{Z}}}

\nc{\sA}{{\mathsf{A}}}
\nc{\sB}{{\mathsf{B}}}
\nc{\sC}{{\mathsf{C}}}
\nc{\sD}{{\mathsf{D}}}
\nc{\sF}{{\mathsf{F}}}
\nc{\sK}{{\mathsf{K}}}
\nc{\sM}{{\mathsf{M}}}
\nc{\sO}{{\mathsf{O}}}
\nc{\sQ}{{\mathsf{Q}}}
\nc{\sP}{{\mathsf{P}}}
\nc{\sZ}{{\mathsf{Z}}}
\nc{\sfp}{{\mathsf{p}}}
\nc{\sr}{{\mathsf{r}}}
\nc{\sfb}{{\mathsf{b}}}
\nc{\sfc}{{\mathsf{c}}}
\nc{\sd}{{\mathsf{d}}}

\nc{\tA}{{\widetilde{\mathbf{A}}}}
\nc{\tB}{{\widetilde{\mathcal{B}}}}
\nc{\tg}{{\widetilde{\mathfrak{g}}}}
\nc{\tG}{{\widetilde{G}}}
\nc{\TM}{{\widetilde{\mathbb{M}}}{}}
\nc{\tO}{{\widetilde{\mathsf{O}}}{}}
\nc{\tU}{{\widetilde{\mathfrak{U}}}{}}
\nc{\TZ}{{\tilde{Z}}}
\nc{\tx}{{\tilde{x}}}
\nc{\tbv}{{\tilde{\bv}}}
\nc{\tfP}{{\widetilde{\mathfrak{P}}}{}}
\nc{\tz}{{\tilde{\zeta}}}
\nc{\tmu}{{\tilde{\mu}}}

\nc{\urho}{\underline{\rho}}
\nc{\uB}{\underline{B}}
\nc{\uC}{{\underline{\mathbb{C}}}}
\nc{\ui}{\underline{i}}
\nc{\uj}{\underline{j}}
\nc{\ofP}{{\overline{\mathfrak{P}}}}
\nc{\oB}{{\overline{\mathcal{B}}}}
\nc{\og}{{\overline{\mathfrak{g}}}}
\nc{\oI}{{\overline{I}}}

\nc{\eps}{\varepsilon}
\nc{\hrho}{{\hat{\rho}}}

\nc{\one}{{\mathbf{1}}}
\nc{\two}{{\mathbf{t}}}

\nc{\Rep}{{\mathop{\operatorname{\rm Rep}}}}
\nc{\Tot}{{\mathop{\operatorname{\rm Tot}}}}
\nc{\Ker}{{\mathop{\operatorname{\rm Ker}}}}
\nc{\Hilb}{{\mathop{\operatorname{\rm Hilb}}}}
\nc{\End}{{\mathop{\operatorname{\rm End}}}}
\nc{\Ext}{{\mathop{\operatorname{\rm Ext}}}}
\nc{\CHom}{{\mathop{\operatorname{{\mathcal{H}}\it om}}}}
\nc{\GL}{{\mathop{\operatorname{\rm GL}}}}
\nc{\gr}{{\mathop{\operatorname{\rm gr}}}}
\nc{\Id}{{\mathop{\operatorname{\rm Id}}}}
\nc{\de}{{\mathop{\operatorname{\rm def}}}}
\nc{\length}{{\mathop{\operatorname{\rm length}}}}
\nc{\supp}{{\mathop{\operatorname{\rm supp}}}}

\nc{\Cliff}{{\mathsf{Cliff}}}
\nc{\Fl}{{\mathsf{Fl}}}
\nc{\Fib}{{\mathsf{Fib}}}
\nc{\Coh}{{\mathsf{Coh}}}
\nc{\FCoh}{{\mathsf{FCoh}}}

\nc{\reg}{{\text{\rm reg}}}

\nc{\cplus}{{\mathbf{C}_+}}
\nc{\cminus}{{\mathbf{C}_-}}
\nc{\cthree}{{\mathbf{C}_*}}
\nc{\Qbar}{{\bar{Q}}}

\nc{\bh}{{\bar{h}}}
\nc{\bOmega}{{\overline{\Omega}}}

\nc{\seq}[1]{\stackrel{#1}{\sim}}

\nc{\bSet}{{\mathbf Set}}
\nc{\BSet}{{\mathbb Set}}
\nc{\BVect}{{\mathbb Vect}}
\nc{\wh}{\widehat}

\begin{document}

\title[Irreducibility of certain induced representations]
{Algebraic groups over a 2-dimensional local field: \\
Irreducibility of certain induced representations}

\author{Dennis Gaitsgory and David Kazhdan}

\dedicatory{To Raoul Bott, with admiration}

\begin{abstract}
Let $G$ be a split reductive group over a local field $\bK$, and let
$G((t))$ be the corresponding loop group. In \cite{GK} we have introduced
the notion of a representation of (the group of $\bK$-points) of $G((t))$
on a pro-vector space. In addition, we have defined an induction
procedure, which produced $G((t))$-representations from usual
smooth representations of $G$. We have conjectured that the induction
of a cuspidal irreducible representation of $G$ is irreducible.
In this paper we prove this conjecture for $G=SL_2$.

\end{abstract}

\maketitle

\section{The result}

\ssec{Notation}

The notation in this paper follows closely that of  \cite{GK}. Let remind
the main characters. We denote by $Set_0$ the category of finite sets,
and $\bSet:=\on{Ind}(\on{Pro}(Set_0))$, $\BSet=\on{Ind}(\on{Pro}(\bSet))$.
By $Vect_0$ we denote the category of finite-dimensional vector spaces
over $\BC$, $Vect=\on{Ind}(Vect_0)$ is the category of all vector spaces,
and $\BVect$ is the category $\on{Pro}(Vect)$ of pro-vector spaces.

Let $G$ be a split reductive group over $\bK$; $\bG$ the
corresponding group-object of $\bSet$. We have the pro-algebraic group of
arcs $G[[t]]$ and for $n\in \BN$ we denote by $G^n\subset G[[t]]$
the corresponding congruence subgroup. By $\bG[[t]]$ (resp.,
$\bG^n\subset \bG[[t]]$) we denote the corresponding group-objects
of $\on{Pro}(\bSet)$.

Finally $\BG=\bG((t))$ is the group-object of $\BSet$, which is our main
object of study. We denote by $\on{Rep}(\BG)$ the category of
representations of $\BG$ on $\BVect$, cf. \cite{GK}, Sect. 2.

\ssec{}

Let us recall the formulation of Conjecture 4.7 of \cite{GK}.
Recall that we have an exact
functor $r^\BG_\bG:\on{Rep}(\BG)\to \on{Rep}(\bG,\BVect)$,
and its right adjoint, denoted $i^\BG_\bG$ and called the
induction functor.

The functors $r^\BG_\bG$ and $i^\BG_\bG$
are direct loop-group analogs of the Jacquet and induction
functors for usual reductive groups over $\bK$.

\medskip

Let $\pi$ be an irreducible cuspidal representation of $\bG$, and
set $\Pi:=i^\BG_\bG(\pi)$. In \cite{GK}, Sect. 4.5 it was shown that
the cuspidality assumption on $\pi$ implies that the natural map
\begin{equation} \label{jacquet of ind}
r^\BG_\bG(\Pi)=r^\BG_\bG\circ i^\BG_\bG(\pi)\to \pi
\end{equation}
is an isomorphism. In particular, this implies that
$$\on{End}_{\on{Rep}(\BG)}(\Pi)\simeq
\on{Hom}_{\on{Rep}(\bG,\BVect)}(r^\BG_\bG(\Pi),\pi)\simeq
\on{Hom}_{\on{Rep}(\bG,\BVect)}(\pi,\pi)\simeq \BC.$$

We have formulated:

\begin{conj}  \label{irr}
The object $\Pi\in \on{Rep}(\BG)$ is irreducible.
\end{conj}

In this paper we will prove:

\begin{thm}  \label{main}
\conjref{irr} holds for $G=SL_2$.
\end{thm}

Note that in \cite{GK} \conjref{irr} was stated slightly more
generally, when we allow representations of a central extension
$\wh\BG$ with a given central charge. The proof of \thmref{main}
generalizes to this set-up in a straightforward way.

\medskip

It should be remarked that from the definition of the category of
representations of $G((t))$, it is not at all clear that $G((t))$ admits any
non-trivial irreducible representations is non-obvious. Therefore,
the fact that the above-mentioned irreducibility conjecture holds
is somewhat surprising.

\ssec{}

We will now consider a functor $\on{Rep}(\bG,\BVect)\to \on{Rep}(\BG)$,
which will be the {\it left} adjoint of the functor $r^\BG_\bG$.

First, recall from \cite{GK1}, Proposition 2.7, that the functor
$\on{Coinv}_{\bG^1}: \on{Rep}(\bG^1,\BVect)\to \BVect$ does admit
a left adjoint, denoted $\on{Inf}^{\bG^1}$.

\begin{prop}  \label{inf}
The functor $\on{Coinv}_{\bG^1}:\on{Rep}(\bG[[t]],\BVect)\to
\on{Rep}(\bG,\BVect)$ admits a left adjoint.
\end{prop}

\begin{proof}

For $\pi=(\BV,\rho)\in \on{Rep}(\bG,\BVect)$, consider
the functor $\on{Rep}(\bG[[t]],\BVect)\to Vect$ given by
$$\Pi\mapsto \on{Hom}_{\on{Rep}(\bG,\BVect)}(\pi,
\on{Coinv}_{\bG^1}(\Pi)).$$

We claim that it is enough to show that this functor is
pro-representable. Indeed, this follows by combining Lemma 1.2,
Proposition 2.5 an Lemma 2.7 of \cite{GK}.

Consider the object $\on{Inf}^{\bG^1}(\BV)\in
\on{Rep}(\bG^1,\BVect)$, where $\BV$ is regarded just as
a pro-vector space, and
$$\on{Coind}^{\bG[[t]]}_{\bG^1}(\on{Inf}^{\bG^1}(\BV))\in
\on{Rep}(\bG[[t]],\BVect),$$
where $\on{Coind}^{\bG[[t]]}_{\bG^1}$ is as in \cite{GK1},
Corollary 2.34.

Evidently,
$$\on{Hom}_{\on{Rep}(\bG,\BVect)}(\pi,
\on{Coinv}_{\bG^1}(\Pi))\hookrightarrow
\on{Hom}_{\BVect}(\BV,\on{Coinv}_{\bG^1}(\Pi)),$$
and the latter, in turn, identifies with
$$\on{Hom}_{\on{Rep}(\bG^1,\BVect)}\left(\on{Inf}^{\bG^1}(\BV),\Pi\right)
\simeq \on{Hom}_{\on{Rep}(\bG[[t]],\BVect)}
\left(\on{Coind}^{\bG[[t]]}_{\bG^1}(\on{Inf}^{\bG^1}(\BV)),\Pi\right).$$

Hence, the pro-representability follows from Proposition 1.4
of \cite{GK}.

\end{proof}

We will denote the resulting functor by
$\on{Inf}^{\bG[[t]]}_{\bG}$. Note that by construction,
for a representation $\pi$ of $\bG$ we have a surjection
$$\on{Coind}^{\bG[[t]]}_{\bG^1}(\on{Inf}^{\bG^1}(\pi))
\twoheadrightarrow \on{Inf}^{\bG[[t]]}_{\bG}(\pi).$$

By composing $\on{Inf}^{\bG[[t]]}_{\bG}$ with the functor
$\on{Coind}^\BG_{\bG[[t]]}:
\on{Rep}(\bG[[t]],\BVect)\to \on{Rep}(\BG)$ we obtain a functor,
left adjoint to $r^\BG_\bG$.

\medskip

We will now formulate the main step in the proof of \thmref{main}.
Note that if $\pi$ is a cuspidal representation of $\bG$,
isomorphism \eqref{jacquet of ind} implies that we have
a canonical map
\begin{equation} \label{key map}
\on{Coind}^\BG_{\bG[[t]]}(\on{Inf}^{\bG[[t]]}_{\bG}(\pi))\to \Pi.
\end{equation}

We will deduce \thmref{main} from the following one:

\begin{thm} \label{surj}
If $G=SL_2$, the map of \eqref{key map} is surjective.
\end{thm}

Of course, we conjecture that the map \eqref{key map} is
surjective for any $G$, but we are unable to prove that
at the moment.

\ssec{}

Let us show how \thmref{surj} implies \thmref{main}.
Suppose that $\Pi'$ is a non-zero sub-object of $\Pi$ and let
$\Pi'':=\Pi/\Pi'$ be the quotient. By definition of the induction
functor, we have a map in $\on{Rep}(\bG,\BVect)$.
$$r^{\BG}_\bG(\Pi')\to \pi.$$

Using Proposition 2.4. of \cite{GK}, we obtain that
$r^{\BG}_\bG(\Pi')$ must surject onto $\pi$, since the latter was
assumed irreducible. Since the functor $r^{\BG}_\bG$ is exact (cf.
Lemma 2.6. of {\it loc.cit.}), and because of isomorphism
\eqref{jacquet of ind}, this implies that $r^{\BG}_\bG(\Pi'')=0$.

However,
$\on{Hom}_{\on{Rep}(\BG)}
(\on{Coind}^{\BG}_{\bG[[t]]}(\on{Inf}^{\bG[[t]]}_\bG(\pi), \Pi'')\simeq
\on{Hom}_{\on{Rep}(\bG,\BVect)}(\pi,r^{\BG}_\bG(\Pi'')$. By \thmref{surj},
this implies that $\Pi''=0$.

\section{The key lemma}

\ssec{}

The rest of the paper is devoted to the proof of \thmref{surj}.
We will slightly abuse the notation, and for a scheme $Y$ over $\bK$
we will make no distinction between the corresponding object
$\bY\in \bSet$ and $Y(\bK)$, regarded as a topological space.

Recall the affine Grassmannian
$\on{Gr}_G=G((t))/G[[t]]$ of $G$, and the corresponding object
${\mathbf {Gr}}_G\in \on{Ind}(\bSet)$. Let us represent $\on{Gr}_G$ as the
direct limit of closures of $G[[t]]$-orbits, $\ol{\Gr}{}_G^\lambda$, with
respect to the natural partial ordering on the set of dominant coweights.

Let us also denote by $\wt{\Gr}_G$ the ind-scheme $G((t))/G^1$, which
is a principal $G$-bundle over $\Gr_G$. Let $\wt\Gr{}^\lambda_G$ and
$\wt{\ol\Gr}{}^\lambda_G$ denote the preimages in $\wt\Gr_G$ of the
$G[[t]]$-orbit $\Gr^\lambda_G$ and its closure, respectively. Let
$\wt{\mathbf {Gr}}{}^\lambda_G$ and $\wt{\ol{\mathbf {Gr}}}{}^\lambda_G$
denote the corresponding objects of $\bSet$.

\medskip

By construction (cf. \cite{GK}, Sect. 3.9), as a $\bG[[t]]$-representation, $\Pi$ is
the inverse limit of $\ol{\Pi}^\lambda$, where each $\ol{\Pi}^\lambda$
is the vector space consisting of locally constant $\bG$-equivariant functions
on $\wt{\ol{\mathbf {Gr}}}{}^\lambda_G$ with values in $\pi$.

Set $\Pi^\lambda$ be the kernel of
$\ol{\Pi}^\lambda\to \underset{\lambda'<\lambda}\oplus\, \ol{\Pi}^{\lambda'}$.
Let $ev$ denote the natural evaluation map $\Pi\to \Pi^0\simeq \pi$,
which sends a function $f\in \on{Funct}^{lc}(\wt{\ol{\mathbf {Gr}}}{}^\lambda_G,\pi)$
to $f(1)$. More generally, for $\wt{g}\in \wt{\ol{\mathbf {Gr}}}{}^\lambda_G$,
we will denote by $ev_{\wt{g}}$ the map $\Pi\to \pi$, corresponding to evaluation
at $\wt{g}$.

\medskip

To prove \thmref{surj} we must show that the composition
\begin{equation}  \label{surj lambda}
\on{Coind}^{\BG}_{\bG[[t]]}\left(\on{Inf}^{\bG[[t]]}_\bG(\pi)\right)\to \Pi\to
\ol{\Pi}^{\lambda}
\end{equation}
is surjective for every $\lambda$. We will argue by induction.
Therefore, let us first check that the map of \eqref{surj lambda}
is indeed surjective for $\lambda=0$.

We have a natural map
\begin{equation} \label{infl to coind}
\on{Inf}^{\bG^1}(\pi)\to \on{Inf}^{\bG[[t]]}_{\bG}(\pi)\to
\on{Coind}^{\BG}_{\bG[[t]]}\left(\on{Inf}^{\bG[[t]]}_\bG(\pi)\right),
\end{equation}
and its composition with
$$\on{Coind}^{\BG}_{\bG[[t]]}\left(\on{Inf}^{\bG[[t]]}_\bG(\pi)\right)\to
\Pi\overset{\on{ev}}\to\pi$$
is the natural surjection $\on{Inf}^{\bG^1}(\pi)\to \pi$.

\medskip

Thus, we have to carry out the induction step. We
will suppose that the composition
$$\on{Coind}^{\BG}_{\bG[[t]]}\left(\on{Inf}^{\bG[[t]]}_\bG(\pi)\right)\to \Pi\to
\ol{\Pi}^{\lambda'}$$
is surjective for $\lambda'<\lambda$, and we must show that
\begin{equation} \label{int surj}
\on{Coind}^{\BG}_{\bG[[t]]}\left(\on{Inf}^{\bG[[t]]}_\bG(\pi)\right)
\underset{\ol{\Pi}^\lambda}\times \Pi^\lambda\to \Pi^\lambda
\end{equation}
is surjective as well.

\ssec{}

For $\lambda$ as above let $t^\lambda$ be the corresponding
point in $\bG((t))$. By a slight abuse of notation we will denote by the
same symbol its image in $\Gr_G$ and $\wt\Gr_G$.

Consider the action of $G^1\subset G((t))$
on $\Gr_G$ given by
$$g\times x=\on{Ad}_{t^\lambda}(g)\cdot x.$$
Let $Y\subset \Gr_G$ be the closure of
$\on{Ad}_{t^\lambda}(G^1)\cdot \ol{\Gr}{}_G^\lambda$.
Let $G_\lambda$ be a finite-dimensional quotient of $G^1$,
through which it acts on $Y$.

We will denote by $\bY$ and $\bG_\lambda$,
respectively, the corresponding objects of $\bSet$.
Let $\Pi_Y$ denote the quotient of $\Pi$, equal to the space
of $\bG$-equivariant locally constant $\pi$-valued functions
on the set of $\bK$-points of the preimage of $Y$ in $\wt\Gr_G$.

\medskip

Let $N\subset G$ be the maximal unipotent subgroup.
Since $\lambda$ is dominant, $\on{Ad}_{t^\lambda}(N[[t]])$
is a subgroup of $N[[t]]$. Let $N^\lambda\subset N[[t]]$
be any normal subgroup of finite codimension, contained in
$\on{Ad}_{t^\lambda}(N[[t]])$.
(Later we will specify to the case when $G=SL_2$; then $N\simeq G_a$
and is abelian, and we will take $N^\lambda=\on{Ad}_{t^\lambda}(N[[t]])$.)
Let $N_\lambda$ denote the quotient $N[[t]]/N^\lambda$, and let
$\bN_\lambda$ be the corresponding group-object in $\bSet$.

Let now $K_{\bN}$ be an open compact subgroup in
$\bN_\lambda$, and $K_{\bG_\lambda}$ an open compact subgroup
in $\bG_\lambda$.

\medskip

Now we are ready to state our main technical claim, \mainlemref{integral}.
However, before doing that, let us explain the idea behind this lemma:

From the isomorphism \eqref{jacquet of ind}, we will obtain that for
any $f\in \Pi_Y$ and a {\it large} compact subgroup $K_{\bG_\lambda}$
as above, the integral $f':=\underset{k\in K_{\bG_\lambda}}\int f^k$ "localizes"
near $t^\lambda$, i.e., $f'$ will be $0$
outside a "small" ball around $t^\lambda$. We will then average $f'$ with
respect to a fixed open subgroup $K_{\bN}$ of $\bN_\lambda$, and obtain a
new element, denoted $f''\in \ol\Pi{}^\lambda$.

\mainlemref{integral} will insure that the compact subgroup $K_{\bG_\lambda}$
can be chosen so that $f''$ will still be localized near $t^\lambda$, and such the
resulting elements $f''$ for various subgroups $K_{\bN}$, and their translations
by elements of $\bG((t))$, span $\Pi^\lambda$.

\ssec{}

In precise terms, we proceed as follows.
Consider the operator $A_{K_\bN,K_{\bG_\lambda}}:\Pi\to \pi$ given by
$$f\mapsto \underset{n\in K_\bN}\int \underset{k\in K_{\bG_\lambda}}\int\,
ev_{t^\lambda}(f^{n\cdot k}),$$
where the integral is taken with respect to the Haar measures on both groups.
(In the above formula $f\mapsto f^x$ denotes the action of $x\in \bG((t))$ on $\Pi$.)
By the definition of $\Pi_Y$, the above map factors through $\Pi\twoheadrightarrow \Pi_Y\to \pi$.

\medskip

For a point $\wt{g}\in \wt{\ol{\mathbf {Gr}}}{}_G^\lambda$ we have a map
$A_{\wt{g},K_{\bG_\lambda}}:\Pi\to \pi$ given by
$$f\mapsto \underset{k\in K_{\bG_\lambda}}\int\, ev_{\wt{g}}(f^{k}).$$
This map also factors through $\Pi_Y$.

\medskip

Our main technical claim, which we prove for
$G=SL_2$ is the following.
(We do not know whether an analogous statement holds for
groups $G$ of higher rank.)

\begin{mainlemma}  \label{integral}
For $v\in \pi$, an open compact subgroup
$K_\bN\subset \bN_\lambda$
and open compact subset $\bX\subset {\mathbf {Gr}}_G^\lambda$
containing $t^\lambda$, there exists a finite-dimensional
subspace $\sF(v)\subset \Pi_Y$ and an increasing exhausting
family of compact subgroups $K^\alpha_{\bG_\lambda}(v)\subset \bG_\lambda$
such that:

\smallskip

\noindent{\em(1)} For all sufficiently large indices $\alpha$
the vector $v$ would
belong to the image of $A_{K_\bN,K_{\bG_\lambda}^\alpha(v)}(\sF(v))$.

\smallskip

\noindent{\em(2)} For every $f\in \sF(v)$ and for all
sufficiently large indices $\alpha$, the vector
$A_{\wt{g},K_{\bG_\lambda}^\alpha(v)}(f)$ will vanish, unless the image of
$\wt{g}$ under $ \wt{\ol{\mathbf {Gr}}}{}_G^\lambda\to  \ol{\mathbf {Gr}}{}_G^\lambda$
belongs to $\bX$.

\end{mainlemma}

\ssec{}

Let us show how \lemref{integral} implies the induction step
in the proof of \thmref{surj}.

\medskip

Recall that the orbit of the point $t^\lambda$ under the action of $N[[t]]$
is open in $\Gr_G^\lambda$. For an open compact subgroup $K_\bN\subset
\bN_\lambda$, let $\bX\subset {\mathbf {Gr}}^\lambda_G$
be its orbit under $K_\bN$. Let $(\ol{\Pi}^\lambda)^{K_\bN}\subset
\ol{\Pi}^\lambda$ be the subspace of $K_\bN$-invariants. We have a direct sum
decomposition
$$(\ol{\Pi}^\lambda)^{K_\bN}=\bV_1\oplus \bV_2,$$
where the first direct summand consists of functions that vanish on the
preimage of $\bX$,
and the second one functions of sections that vanish outside the preimage
of $\bX$. We have $\bV_2\subset \Pi^\lambda$ and
the map $ev_{t^\lambda}$ identifies $\bV_2$ with $\pi$.

\medskip

We claim that it suffices to show that the image of the map
\begin{equation} \label{map with averaging}
\on{Coind}^{\BG}_{\bG[[t]]}\left(\on{Inf}^{\bG[[t]]}_\bG(\pi)\right)
\to\Pi\to\ol{\Pi}^\lambda\to (\ol{\Pi}^\lambda)^{K_\bN},
\end{equation}
where the last arrow is given by averaging with respect to $K_\bN$, contains
$\bV_2$.

Indeed, let $G[[t]]_\lambda$ be a finite-dimensional quotient through which
$G[[t]]$ acts on $\Gr_G^\lambda$, and let $\bG[[t]]_\lambda$ be the corresponding
group-object of $\bSet$. The vector space $\Pi^\lambda$ is spanned by elements
of the following form.
Each is invariant under some (small) open compact subgroup
$K_{\bG[[t]]_\lambda}\subset \bG[[t]]_\lambda$, and is supported on a preimage in
$\wt{\mathbf {Gr}}{}^\lambda_G$ of a single $K_{\bG[[t]]_\lambda}$-orbit on
${\mathbf {Gr}}^\lambda_G$.  By $\bG[[t]]$-invariance, we can assume that
the orbit in question is that of the element $t^\lambda\in {\mathbf {Gr}}^\lambda_G$.

By setting $K_{\bN}:=\bN_\lambda\cap K_{\bG[[t]]_\lambda}$, we obtain that any
element of the form specified above is contained in the corresponding $\bV_2$.

\medskip

We will show that \mainlemref{integral} implies that
$\bV_2$ belongs to the image of the map
$$\on{Inf}^{\bG^1}(\pi)\to
\on{Coind}^{\BG}_{\bG[[t]]}\left(\on{Inf}^{\bG[[t]]}_\bG(\pi)\right)\to
(\ol{\Pi}^\lambda)^{K_\bN},$$
where first the arrow is the composition of the map of \eqref{infl to coind},
followed by the action of $t^\lambda$.

For that let us write down in explicit terms the composition
\begin{equation}  \label{map from inflation}
\on{Inf}^{\bG^1}(\pi)\to
\on{Coind}^{\BG}_{\bG[[t]]}\left(\on{Inf}^{\bG[[t]]}_\bG(\pi)\right)\to
\Pi\to \Pi_Y.
\end{equation}

First, let us observe that the resulting map factors through
the surjection $\on{Inf}^{\bG^1}(\pi)\twoheadrightarrow
\on{Inf}^{\bG_\lambda}(\pi)$.
Secondly, let us recall (cf. \cite{GK1}, Sect. 2.8) that
$\on{Inf}^{\bG_\lambda}(\pi)$ is the inductive limit, taken in $\BVect$,
over finite-dimensional subspaces $\sF'\subset \pi$ of
$$\underset{\alpha}{\underset{\longleftarrow}{"lim"}}\,
\on{Distr}_c(\bG_\lambda/K_{\bG_\lambda}^\alpha)\otimes \sF',$$
where $K_{\bG_\lambda}^\alpha$ runs through any exhausting family of
open compact subgroups of $\bG_\lambda$.

By \eqref{jacquet of ind},
the map $\Pi_Y\overset{ev_{t^\lambda}}\to\pi$ induces an isomorphism
$\on{Coinv}_{\bG_\lambda}(\Pi_Y)\simeq \pi$.
For a given finite-dimensional subspace $\sF'$ let us choose
a finite-dimensional subspace $\sF\subset \Pi_Y$ which projects
surjectively onto $\sF'$, and for every index $\alpha$ consider the map
$$\on{Distr}_c(\bG_\lambda/K_{\bG_\lambda}^\alpha)\otimes \sF\to\Pi_Y$$
given by $$\mu\otimes f\mapsto \mu\star f,$$
where $f\in \sF$ and
$\mu\in \on{Distr}_c(\bG_\lambda/K_{\bG_\lambda}^\alpha)$ is regarded
as an element of the Hecke algebra.

The resulting system of maps (eventually in $\alpha$) factors through
$\on{Distr}_c(\bG_\lambda/K_{\bG_\lambda}^\alpha)\otimes \sF\twoheadrightarrow
\on{Distr}_c(\bG_\lambda/K_{\bG_\lambda}^\alpha)\otimes \sF'$, and
defines the map in \eqref{map from inflation}.

\medskip

Let us now recall that if $\BW=\underset{\longleftarrow}{lim}\, \bW_\alpha$
is a pro-vector
space mapping to a vector space $\bV$, the surjectivity of this map means
that the eventually defined maps $\bW_\alpha\to \bV$ are all surjective, or,
which is the same, that $\forall v\in \bV$, $v\in \on{Im}(\bW_\alpha)$
for those indices $\alpha$, for which the map $\BW\to \bV$ factors through
$\bW_\alpha\to \bV$.

For a vector $v\in \pi$, let $\sF(v)$ be the finite-dimensional
subspace of $\Pi_Y$, given by \lemref{integral},
and let $K^\alpha_{\bG_\lambda}(v)$ be the corresponding system of subgroups.
Let $\sF'(v)$ denote the image of $\sF(v)$ in $\pi$.

Consider the composition:
$$\on{Distr}_c(\bG_\lambda/K_{\bG_\lambda}^\alpha(v))\otimes \sF(v)\to
\Pi_Y\to \ol{\Pi}^\lambda\to (\ol{\Pi}^\lambda)^{K_\bN}.$$

Let us take the unit element in
$\on{Distr}_c(\bG_\lambda/K_{\bG_\lambda}^\alpha(v))$, corresponding
to the Haar measure on $K_{\bG_\lambda}^\alpha(v)$. We obtain a map
$\sF(v)\to (\ol{\Pi}^\lambda)^{K_\bN}$.

By \lemref{integral}(2), the image of this map is contained
in $\bV_2$. When we further compose it with the evaluation map
$\bV_2\hookrightarrow \Pi_Y\to \pi$ we obtain a map
$\sF(v)\to \pi$ equal to $A_{K_\bN,K_{\bG_\lambda}^\alpha(v)}$,
whose image contains $v$, by \lemref{integral}(1).

This establishes the required surjectivity.

\section{Proof of \mainlemref{integral}}

\ssec{}

For a given subgroup $K_\bN\subset
\bN_\lambda$, a subset
$\bX\subset {\mathbf {Gr}}_G^\lambda$ and an arbitrary
finite-dimensional subspace $\sF\subset \Pi_Y$ we will
construct a family of open compact subgroups
$K_{\bG_\lambda}\subset \bG_\lambda$, such that
the expressions
$$A_{K_\bN,K_{\bG_\lambda}}(f) \text{ and }
A_{\wt{g},K_{\bG_\lambda}}(f)$$
for $f\in \sF$ can be evaluated explicitly.

\medskip

From now on we will fix $G=SL_2$. We will  change
the notation slightly, and identify the set
of dominant coweights with $\BN$; in which case we will
replace $\lambda$ by $l$ and $t^\lambda\in \bG((t))$ becomes
the matrix
$$\begin{pmatrix} t^l & 0 \\
0 & t^{-l}
\end{pmatrix}.
$$

Let us translate our initial subscheme $Y$ by $t^{-\lambda}$,
in which case the point $t^\lambda$ itself will go over to
the unit point $1_{\Gr_G}\in \Gr_G$, and $t^{-\lambda}\cdot Y$ will
be contained in $\ol\Gr{}_G^{2l}$. (We denote by $\Gr_G^r$ the
$G[[t]]$-orbit of the point
$\begin{pmatrix} t^r & 0 \\
0 & t^{-r}
\end{pmatrix}$ in $\Gr_G$, and by $\ol\Gr{}_G^r$ its closure.)
For the purposes of \lemref{integral}
we can replace $t^{-\lambda}\cdot Y$
by the entire $\ol\Gr{}_G^{2l}$, with the standard
action of the congruence subgroup $G^1$.

Note that the action of
$G^1$ on $\ol\Gr{}_G^{2l}$ (resp., $\wt{\ol\Gr}{}_G^{2l}$)
factors through $G^1/G^{4l}$ (resp., $G^1/G^{4l+1}$).

\medskip

For an integer $r$ let us denote by $G_r$ the quotient
$G^1/G^{2r+1}$, and by $N_r$ the quotient $t^{-r}\cdot N[[t]]/N[[t]]$.
We will write elements of $\bN_r$ as
$\underset{1\leq i\leq r}\Sigma\, t^{-i}\cdot n_i$ with $n_i\in \bK$,  and thus
think of it as an $r$-dimensional vector space over $\bK$.

Similarly, we will identify $\bG_r:=\bG^1/\bG^{2r+1}$ with an $6r$-dimensional
vector space over $\bK$, by writing its elements as matrices:
$$
\on{Id}+
\begin{pmatrix}
k_{11} & k_{12} \\
k_{21} & k_{22}
\end{pmatrix}
$$
and $k_{lm}=\underset{1\leq i\leq 2r}\Sigma\,t^{i}\cdot (k_{lm})_i$. In particular,
we can speak of $O_\bK$-lattices in $\bG_r$, where $O_\bK\subset \bK$
is the ring of integers.

\ssec{}

In what follows, for a point $g\in \bG((t))$, we will denote by $\wt{g}$
(resp., $\ol{g}$) its image in $\wt{\mathbf {Gr}}_G$ (resp., ${\mathbf {Gr}}_G$).

Thus, we are interested in computing the integral
$$\underset{k\in K_{\bG_r}}\int\, ev(f^{g\cdot k}),$$
when $g$ is such that
either $g\in \bN_r$, or the corresponding point
$\ol{g}\in {\mathbf {Gr}}_G$ lies in $\ol{\mathbf {Gr}}{}_G^r-\bX$,
where $\bX$ is a fixed open compact subset of $\ol{\mathbf {Gr}}{}_G^r$
containing $1_{\Gr_G}$, and $f\in \sF$, where $\sF$ is a fixed
finite-dimensional subspace of $\ol\Pi^r$.

\medskip

Let $\fp$ denote the projection $\wt\Gr_G\to \Gr_G$.
Let $\fs$ be a continuous section $\ol{\mathbf {Gr}}{}^r_G\to
\wt{\ol{\mathbf {Gr}}}{}^r_G$, such that $\fs(1_{\Gr_G})=1_{\wt{\Gr}_G}$.
A choice of such section defines an isomorphism
$$\wt{\ol{\mathbf {Gr}}}{}^r_G\simeq \ol{\mathbf {Gr}}{}^r_G\times \bG.$$
We will denote by $\fq$ the resulting map $\wt{\ol{\mathbf {Gr}}}{}^r_G\to \bG$.

Let us fix an open neighbourhood $\bZ$ of $1_{\Gr_G}$ in $\ol{\mathbf {Gr}}{}_G^r$
small enough so that $$f(\fs(x))=ev(f)$$ for $x\in \bZ$ and
$f\in \sF$. Let $K_\bG(\sF)$ be an open compact subgroup of $\bG$,
such that $ev(f)\in \pi$ is $K_{\bG}(\sF)$-invariant for $f\in \sF$.

Let $K_{\bN_r}$ be an open compact subgroup of $\bN_r$.

\begin{prop}   \label{construction of sgr}
There exists an $O_\bK$-lattice
$K_{\bG_r}\subset \bG_r$, which contains any given open
subgroup of $\bG_r$, such that the following is satisfied:

\smallskip

\noindent (1) There exists an open compact subgroup
$K_{\bN_r}^{sm}\subset K_{\bN_r}$ such that:

\begin{itemize}

\item(1a) For $g=k\cdot n\in \bG((t))$ with
$k\in K_{\bG_r}$ and $n\in K_{\bN_r}^{sm}$, the corresponding point
$\ol{g}\in \ol{\mathbf {Gr}}{}_G^{l}$ belongs to $\bZ$.

\item(1b)
For $g$ as above, the left coset of $\fq(\wt{g})\in \bG$
with respect to $K_\bG(\sF)\subset \bG$ equals that of
$$
\begin{pmatrix}
1 & -\underset{1\leq i\leq r}\Sigma\, (k_{12})_{2i}\cdot n_i^2 \\
0 & 1.
\end{pmatrix}
$$

\item(1c)
The integral
$\underset{k\in K_{\bG_r}}\int\, ev(f^{n\cdot k})$ vanishes if
$n\in K_{\bN_r}-K_{\bN_r}^{sm}$ and $f\in \sF$.

\end{itemize}

\smallskip

\noindent (2) If $g\in \bG((t))$, such that
$\ol{g}\in \ol{\mathbf {Gr}}{}_G^r-\bX$, the integral
$\underset{k\in K_{\bG_r}}\int\, ev(f^{g\cdot k})$ vanishes.

\end{prop}

\ssec{}

Let us deduce \mainlemref{integral} from \propref{construction of sgr}.  Given
a vector $v\in \pi$ let us first define the subspace $\sF(v)\in \ol{\Pi}^r$.

Recall that $\bN\simeq \bK$ is the maximal unipotent subgroup of $\bG=SL_2(\bK)$,
and let $\bN^*$ denote the Pontriagin dual group. Since $\bN^*$ is also
(non-canonically) isomorphic to $\bK$, we have a valuation map $\nu:\bN^*\to \BZ$,
defined up to a shift. In particular, we can consider the subalgebra
$\on{Funct}_{val}(\bN^*):\simeq \on{Funct}(\BZ)$ inside the algebra $\on{Funct}_{lc}(\bN^*)$
of all locally constant functions on $\bN^*$.

Any smooth representation of $\bN$,
and in particular $\pi$, can be thought of as a module over the algebra of
$\on{Funct}_{lc}(\bN^*)$, such that every
element of this module has compact support.
If a representation is cuspidal, this means that the support of every
section is disjoint from $0\in \bN^*$.

Therefore, if $v$ is a vector in a
cuspidal representation $\pi$, the vector space
$\on{Funct}_{val}(\bN^*)\cdot v\subset \pi$
is finite-dimensional. We denote this vector subspace by $\sF'(v)$ and let
$\sF(v)\subset \ol\Pi^r$ to be any subspace surjecting
onto $\sF'(v)$ by means of $ev$. We claim that $\sF(v)$ satisfies
the requirements of \mainlemref{integral}.

\medskip

Property (2) in the lemma is satisfied due to \propref{construction of sgr}(2).
To check property (1) we will rewrite $A_{K_\bN,K_{\bG_\lambda}}(f)$
more explicitly in terms of the action of $\bG$ on $\pi$.

Note that by \propref{construction of sgr}(1c), the integral
$$\underset{n\in K_{\bN_r}}\int
\underset{k\in K_{\bG_r}}\int\, ev(f^{n\cdot k})$$
equals the integral over a smaller domain, namely,
$$\underset{n\in K_{\bN_r}^{sm}}\int
\underset{k\in K_{\bG_r}}\int\, ev(f^{n\cdot k}).$$

By \propref{construction of sgr}(1a) and (1b), the latter can be rewritten as
\begin{equation} \label{int1}
\underset{n\in K_{\bN_r}^{sm}}\int
\underset{k\in K_{\bG_r}}\int\,
\begin{pmatrix}
1 & \underset{1\leq i\leq l}\Sigma\, (k_{12})_{2i}\cdot n_i^2 \\
0 & 1
\end{pmatrix}
\cdot ev(f).
\end{equation}

For $n=\Sigma\, t^{-i}\cdot n_i\in \bN_r$
consider the map $\phi_n:\bG_r\to \bN$, given by
$$k=\begin{pmatrix}
1+\underset{i}\Sigma\, t^i\cdot (k_{11})_i & \underset{i}\Sigma\, t^i\cdot (k_{12})_i \\
\underset{i}\Sigma\, t^i\cdot (k_{21})_i & 1+\underset{i}\Sigma\, t^i\cdot (k_{22})_i
\end{pmatrix}\mapsto \begin{pmatrix}
1 & \underset{1\leq i\leq r}\Sigma\, (k_{12})_{2i}\cdot n_i^2 \\
0 & 1.
\end{pmatrix}
$$
Thus, the expression in \eqref{int1} can be rewritten as
\begin{equation}  \label{int2}
\underset{n\in K^{sm}_{\bN_r}}\int\,
(\phi_n)_*\left(\mu(K_{\bG_r})\right)\cdot ev(f),
\end{equation}
where $\mu(K_{\bG_r})$ denotes the Haar measure of $K_{\bG_r}$, and
$(\phi_n)_*\left(\mu(K_{\bG_r})\right)$ is its push-forward under $\phi_n$,
regarded as a distribution on $\bN$.

\medskip

Note, however, that when we identify $\bG_r\simeq \bG^1/\bG^{2r+1}$ with a linear space
over $\bK$, the Haar measure on this group goes over to a linear Haar measure.
From this we obtain that for each $n\in \bN_r$, the distribution $(\phi_n)_*\left(\mu(K_{\bG_r})\right)$,
thought of as a function on $\bN^*$, is the characteristic function of some $O_\bK$-lattice in $\bN^*$.
Moreover, this lattice grows as $n\to 0$.

In particular, $(\phi_n)_*\left(\mu(K_{\bG_r})\right)$, as a function on $\bN^*$, belongs
to $\on{Funct}_{val}(\bN^*)$, and the integral $\underset{n\in K^{sm}_{\bN_r}}\int\,
(\phi_n)_*\left(\mu(K_{\bG_r})\right)$, being {\it positive} at every point of $\bN^*$,
defines an invertible element of $\on{Funct}_{val}(\bN^*)$.

Hence $$v\in (\phi_n)_*\left(\mu(K_{\bG_r})\right)\cdot (\on{Funct}_{val}(\bN^*)\cdot v)=
\on{Im} \left(\underset{n\in K^{sm}_{\bN_r}}\int\,
(\phi_n)_*\left(\mu(K_{\bG_r})\right)\cdot ev(\sF(v))\right),$$
which is what we had to show.

\section{Proof of \propref{construction of sgr}}

\ssec{}

We will construct the subgroup $K_{\bG_r}$ by induction with respect
to the parameter $r$. (For property (2) we take $\bX\cap \ol{\mathbf
{Gr}}_G^{r-1}$ as the corresponding open compact subset of
$\ol{\mathbf {Gr}}_G^{r-1}$.)

When $r=0$ all the subgroups in question are trivial. So, we can
assume having constructed the subgroups $K_{\bG_{r-1}}$ and
$K_{\bN_{r-1}}^{sm}$, and let us perform the induction step.
The key observation is provided by the following lemma:

\begin{lem} \label{vanishing}
Let $\bX$ be a compact subset of ${\mathbf Gr}_G^r$
and $f\in \ol\Pi^r$.
Then the integral $\underset{k\in K_{\bG^{2r}/\bG^{2r+1}}}\int\, ev_{\wt{g}}(f^{k})=0$
if $K_{\bG^{2r}/\bG^{2r+1}}$ is a sufficiently large compact subgroup of
$\bG^{2r}/\bG^{2r+1}$, and $\fp(\wt{g})\in \bX$.
\end{lem}

\begin{proof}

Since $\Gr^r_G$ is a $G[[t]]$-orbit of
$t^\lambda:=\begin{pmatrix} t^r & 0 \\
0 & t^{-r}
\end{pmatrix}$
and since $G^{2r}$ is normalized by $G[[t]]$ and acts trivially
on $\ol{\Gr}{}^r_G$, by the compactness of $\bX$,
the assertion of the lemma reduces to the fact that
$$\underset{k\in K_{\bG^{2r}/\bG^{2r+1}}}\int\, ev_{t^\lambda}(f^{k})=0$$ for
every sufficiently large subgroup $K_{\bG^{2r}/\bG^{2r+1}}$.

Note that for $k\in \bG^{2r}$ written as
$$
\begin{pmatrix}
1+t^{2r}\cdot k_{11} & t^{2r}\cdot k_{12} \\
t^{2r}\cdot k_{21} & 1+t^{2r}\cdot k_{22}
\end{pmatrix},
$$
$\on{Ad}_{t^\lambda}(k)\in \bG[[t]]$ projects to the element
$$
\begin{pmatrix}
1 & k_{12} \\
0 & 1
\end{pmatrix} \in \bG=\bG[[t]]/\bG^1.$$
We have
$$ev_{t^\lambda}(f^{k})=\on{Ad}_{t^\lambda}(k)\cdot ev_{t^\lambda}(f).$$

Therefore, the integral in question equals the averaging of the
vector $ev_{t^\lambda}(f)\in \pi$ over a compact subgroup of the
maximal unipotent subgroup of $G$. Moreover, this subgroup grows
together with $K_{\bG^{2r}/\bG^{2r+1}}$. Hence, our assertion
follows from the cuspidality of $\pi$.

\end{proof}

\ssec{}

To carry out the induction step we first choose $K'_{\bG_r}\subset \bG_r$
to be any $O_\bK$-lattice, which projects onto $K_{\bG_{r-1}}
\subset \bG_{r-1}$.

By continuity and the compactness of $K'_{\bG_r}$ there exists
an $O_\bK$-lattice $\bL\subset \bK$, such that for
$K^{sm}_{\bN_r}=K^{sm}_{\bN_{r-1}}+t^{-r}\bL$ the following holds:

\begin{itemize}

\item(a')
For $g=k'\cdot n\in \bG((t))/\bG^1$ with $k\in K'_{\bG_r}$,
$n\in K^{sm}_{\bN_r}$, the point
$\ol{g}\in \ol{\mathbf {Gr}}^r_G$ belongs to $\bZ$.

\item(b')
For $g$ as above, the left coset of $\fq(\wt{g})\in \bG$
with respect to $K_\bG(\sF)\subset \bG$ equals that of
$$
\begin{pmatrix}
1 & -\underset{1\leq i\leq r}\Sigma\, (k_{12})_{2i}\cdot n_i^2 \\
0 & 1.
\end{pmatrix}
$$

\item(c')
For $f\in \sF$, $f(k'\cdot (n'+t^{-r}n_r))=f(k'\cdot n')$ for
$n'\in K^{sm}_{\bN_{r-1}}$, $k'\in K'_{\bG_r}$, and $n_r\in \bL$.

\end{itemize}

Note that for any $n\in \bN_r$, $k'\in \bG_r$ and
$$
k=
\begin{pmatrix}
1+t^{2r}\cdot k_{11} & t^{2r}\cdot k_{12} \\
t^{2r}\cdot k_{21} & 1+t^{2r}\cdot k_{22}
\end{pmatrix}\in \bG^{2r},$$
we have:
\begin{equation} \label{conj}
k\cdot k'\cdot n=k'\cdot n\cdot \begin{pmatrix}
1 & -k_{12}\cdot n_r^2 \\
0 & 1
\end{pmatrix} \,\on{mod} \bG^1.
\end{equation}

\medskip

The group $K_{\bG_r}$ will be obtained from
$K'_{\bG_r}$ by adding to it an (arbitrarily large) lattice
in $\bG^{2r}/\bG^{2r+1}$.

Note that since $G^{2r}$ acts trivially on $\ol\Gr{}^r_G$,
any such subgroup would satisfy condition (1a) of
\propref{construction of sgr}, because
$K_{\bG_r}$ satisfies (a') above. It will also automatically
satisfy (1b) in view of \eqref{conj} and (b') above.
Thus, we have to arrange so that $K_{\bG_r}$
satisfies conditions (1c) and (2) of \propref{construction of sgr}.

\ssec{}

By \lemref{vanishing}, we can find an open compact subgroup
$K_{\bG^{2r}/\bG^{2r+1}}\subset
\bG^{2r}/\bG^{2r+1}$, such that the integrals
$$\underset{k\in K_{\bG^{2r}/\bG^{2r+1}}}\int\,
ev(f^{n\cdot k'\cdot k})$$
would vanish for $f\in \sF$, $k'\in K'_{\bG_r}$ and
$n\in K_{\bN_r}$ is such that $n_r\notin  \bL$.

Let us enlarge the initial $K'_{\bG_r}$ by adding to it
any $O_\bK$-lattice in
$\bG^{2r}/\bG^{2r+1}$ containing the above
$K_{\bG^{2r}/\bG^{2r+1}}$. We claim that the resulting
subgroup satisfies condition (1c) of \propref{construction of sgr}.

\medskip

Indeed, let $n=n'+t^{-r}n_r$, $n'\in \bN_{r-1}, n_r\in \bK$ be
an element in $K_{\bN_r}-K_{\bN_r}^{sm}$. If $n_r\notin \bL$,
the integral vanishes by the choice of $K_{\bG^{2r}/\bG^{2r+1}}$.
Thus, we can assume that $n_r\in \bL$, but
$n'\notin K_{\bN_{r-1}}^{sm}$. But then the integral vanishes
by (c') and the induction hypothesis..

\medskip

Now, let us deal with condition (2) of \propref{construction of sgr}. By the induction
hypothesis, the integrals
\begin{equation} \label{van int}
\underset{k\in K'_{\bG_r}}\int\, ev(f^{g\cdot k})
\end{equation}
vanish when
$\ol{g}\in \ol{\mathbf {Gr}}{}_G^{r-1}-(\ol{\mathbf {Gr}}{}_G^{r-1}\cap \bX)$.

Hence, by continuity and since $K'_{\bG_r}$ is compact,
there exists a neighbourhood $\bX_1$
of $\ol{\mathbf {Gr}}{}_G^{r-1}-(\ol{\mathbf {Gr}}{}_G^{r-1}\cap \bX)$ in
$\ol{\mathbf {Gr}}{}_G^r-\bX$, such that the integral \eqref{van int} will
vanish for the same subgroup $K'_{\bG_r}$ and all $g$ for which
$\ol{g}\in \bX_1$.

The sought-for subgroup $K_{\bG_r}$ will be again obtained from the
initial $K'_{\bG_r}$ by adding to it an arbitrarily large open compact subgroup of
$\bG^{2r}/\bG^{2r+1}$. We claim that for any such $K_{\bG_r}$
the integral
\begin{equation} \label{van int1}
\underset{k\in K_{\bG_r}}\int\, ev(f^{g\cdot k})
\end{equation}
will still vanish for $\ol{g}\in \bX_1$.

This follows from the fact that the $G^{2r}$-action on
$\Gr^r_G$ is trivial, and hence for $k\in \bG^{2r}$,
$k'\in \bG_r$ and $g\in \bG((t))$ projecting to
$\ol{g}\in \ol{\mathbf {Gr}}{}_G^r$,
$$f(k\cdot k'\cdot g)=k_1\cdot f(k'\cdot g)$$
for some $k_1\in \bG$.

\medskip

We choose the suitable subgroup in $\bG^{2r}/\bG^{2r+1}$ as
follows.
Set $\bX_2=(\ol{\mathbf {Gr}}{}_G^r-\bX)-\bX_1$. This is a compact subset
of ${\mathbf {Gr}}{}_G^r$, and let us apply \lemref{vanishing} to the compact set
$K'_{\bG_r}\cdot \bX_2\subset {\mathbf {Gr}}{}_G^r$.

We obtain that there exists an open compact subgroup
$K_{\bG^{2r}/\bG^{2r+1}}\subset \bG^{2r}/\bG^{2r+1}$, such that
$$\underset{k\in K'_{\bG^{2r}/\bG^{2r+1}}}\int\, ev(f^{g\cdot k'\cdot k})=0$$
for $\ol{g}\in \bX_2$, $k'\in K'_{\bG_r}$.

\medskip

Let $K_{\bG_r}$ be the resulting subgroup of $\bG_r$. We claim that
it does satisfy condition (2) of \propref{construction of sgr}.
Indeed, consider again the integral \eqref{van int1} for
$\ol{g}\in \ol{\mathbf {Gr}}{}_G^r-\bX=\bX_1\cup \bX_2$.

We already know that it vanishes for $\ol{g}\in \bX_1$.
And if $\ol{g}\in \bX_2$, it vanishes by the
choice of $K_{\bG_r}$.

This completes the induction step in the proof of 
\propref{construction of sgr}.


\begin{thebibliography}{99}

\bibitem{GK} D.~Gaitsgory and D.~Kazhdan, {\em Representations of
algebraic groups over a 2-dimensional local field}, math.RT/0302174,
to appear in GAFA.

\bibitem{GK1} D.~Gaitsgory and D.~Kazhdan, {\em Algebraic groups over
a 2-dimensional local field: some further constructions}, math.RT/0406282.

\end{thebibliography}
\end{document}